\documentclass[12pt,twoside,reqno]{amsart}
\usepackage{amsxtra,amscd}
\usepackage{graphicx}
\usepackage{amsmath}
\usepackage{amsfonts}
\usepackage{amssymb}

\theoremstyle{definition}

\numberwithin{equation}{section}

\begin{document}
\title{a short note on quasi-galois closed and pseudo-galois covers}
\author{Feng-Wen An}
\address{School of Mathematics and Statistics, Wuhan University, Wuhan,
Hubei 430072, People's Republic of China}
\email{fwan@amss.ac.cn}
\subjclass[2000]{Primary 14J50; Secondary 14G40}
\keywords{arithmetic scheme, automorphism group, Galois cover}
\begin{abstract}
There is a big difference between \textquotedblleft
quasi-galois closed\textquotedblright\ in the eprint (arXiv:0907.0842) and \textquotedblleft
pseudo-galois\textquotedblright\ in the sense of Suslin-Voevodsky. It is nontrivial.
\end{abstract}
\maketitle

Let $X$ and $Y$ be integral schemes of finite types over $Spec(\mathbb{Z}).$
Let $\varphi :X\rightarrow Y$ be a morphism of finite type.

We define quasi-galois
closed covers in \textbf{Definition 1.1} of \cite{An}. If  $\varphi $ is a finite morphism,  Suslin and Voevodsky define pseudo-galois in \textbf{Definition 5.5} of \cite{VS}.

\textbf{Remark.} There is a big difference between \textquotedblleft
quasi-galois closed\textquotedblright\ and \textquotedblleft
pseudo-galois\textquotedblright. This is nontrivial.
For \textquotedblleft quasi-galois closed\textquotedblright , \ the function
field $k\left( X\right) $ is not necessarily algebraic over $k\left(
Y\right) .$ At the same time, for \textquotedblleft
pseudo-galois\textquotedblright , $k\left( X\right) $ must be an algebraic
extension over $k\left( Y\right) .$ 

\textbf{Example.} Let $X=Spec(\mathbb{
Z}[t])$ and $Y=Spec(\mathbb{Z})$ and let $\varphi :X\rightarrow Y$
be the morphism induced by the inclusion. Here $t$ is a variable over $\mathbb{Q}.$
Then $X$ is quasi-galois closed
over $Y.$ It is clear that $X$ is not pseudo-galois over $Y$.

\bigskip

Many thanks for an anonymous referee's comments.

\end{document}